\documentclass[conference]{IEEEtran}
\ifCLASSINFOpdf
\else
\fi
\usepackage{amsthm,amsmath,amssymb}
\usepackage{mathrsfs}
\usepackage{amsfonts}
\usepackage{enumerate}
\usepackage{latexsym}
\usepackage{multirow}
\usepackage{cite}
\usepackage{verbatim}
\usepackage{mathtools}
\usepackage{array}

\usepackage[pdftex]{color,graphicx}

\newcommand{\argmin}{\operatornamewithlimits{\textrm{argmin}}}

\newtheorem{theorem}{Theorem}[section]

\newtheorem{remark}{Remark}[section]

\newcommand{\E}{{\mathbb E}}

\newcommand{\R}{{\mathbb R}}
\renewcommand{\P}{{\mathbb P}}

\newcommand{\F}{{\cal F}}

\newcounter{rcnt}[section]

\def\argmin{\mathop{\rm argmin}}

\begin{document}
%
\title{Information Theoretic validity of Penalized Likelihood}

\author{\IEEEauthorblockN{Sabyasachi Chatterjee}
\IEEEauthorblockA{Department of Statistics\\
Yale University\\
New Haven, Connecticut 06511\\
Email: sabyasachi.chatterjee@yale.edu}
\and
\IEEEauthorblockN{Andrew Barron}
\IEEEauthorblockA{Department of Statistics\\
Yale University\\
New Haven, Connecticut 06511\\
Email: andrew.barron@yale.edu}}


%

\maketitle

\begin{abstract}
Building upon past work, which developed information theoretic notions of when a penalized likelihood procedure can be interpreted as codelengths arising from a two stage code and when the statistical risk of the procedure has a redundancy risk bound, we present new results and risk bounds showing that the $l_1$ penalty in Gaussian Graphical Models fits the above story. We also show how twice the traditional $l_0$ penalty times plus lower order terms which stay bounded on the whole parameter space has a conditional two stage description length interpretation and present risk bounds for this penalized likelihood procedure.
\end{abstract}

\section{Introduction}
It is known that the MDL principle motivates viewing a penalized log likelihood procedure as minimizing the codelengths of a two stage code. Traditionally, this has required that the minimizing space be countable. In past works \cite{BarronHuangLiLuo2008a},\cite{BarronHuangLiLuo2008b} the authors address this issue and  develop a notion as to how to interpret a penalized log likelihood as codelengths arising from a two stage code even when the minimization is done over an uncountable parameter space. We say a given code or its codelengths on a countable sample space are valid if they satisfy the well known Kraft's inequality.
  
In the following $Z$ denotes a sample of data and $\{P_{\theta}: \theta \in \Theta\}$ denotes a family of candidate distributions for $Z.$ The basic idea, showing how codelengths arising from a penalized likelihood procedure, is very simple. Suppose for the given penalty one can construct another valid codelength $L$ on the sample space such that 
\begin{align*}
\min_{\theta \in \Theta} \{-\log P_{\theta}(Z) + pen(\theta)\} \geq L(Z)
\end{align*}
then automatically the penalized likelihood term defines a codelength satisfying Kraft's inequality. We then call such a penalty codelength valid.  The following condition defined in \cite{BarronHuangLiLuo2008a} allows a penalty to be viewed as Kraft summable codelengths:
 Suppose there exists a countable $\tilde{\theta} \subset \Theta$ and a Kraft valid codelength $\L(\tilde{\theta})$ defined on $\tilde{\Theta}$ such that the following holds for any given data $Z,$
\begin{align}{\label{verifycode}}
\begin{split}
&\min_{\theta \in \Theta} \{-\log P_{\theta}(Z) + pen(\theta)\} \geq \\&\min_{\tilde{\theta} \in \tilde{\Theta}} \{-\log P_{\tilde{\theta}}(Z) + \L(\tilde{\theta})\}
\end{split}
\end{align}
then clearly since the right side of~\eqref{verifycode} is a Kraft summable codelength by virtue of having a codelength associated with a two stage code, the left side also is a Kraft summable codelength. It is shown in \cite{BarronHuangLiLuo2008a} and \cite{BarronHuangLiLuo2008b} that the $l_1$ penalty in linear regression and log density estimation problems is indeed codelength valid. Just showing codelength validity of a penalty is  relevant for data compression purposes but not sufficient for generalization guarantees on future data. The authors in \cite{BarronHuangLiLuo2008a} also define a condition for good statistical risk properties to hold. If there exists a countable $\tilde{\theta} \subset \Theta$ and a Kraft valid codelength $\L(\tilde{\theta})$ on it such that the following condition holds:
\begin{align}{\label{verifyrisk}}
\begin{split}
&\min_{\tilde{\theta} \in \F} \left(\log \frac{P_{\tilde{\theta}}(Z)}{P(Z)} - B(P,P_{\theta})  + 2\L(\tilde{\theta}) \right) \leq \\&\min_{\theta \in \Theta} \left(\log \frac{P_{\theta}(Z)}{P(Z)} - B(P,P_{\theta}) + pen(\theta) \right)
\end{split}
\end{align}
then, the right side inherits the positive expectation property from the left side in~\eqref{verifyrisk} as shown in \cite{BarronHuangLiLuo2008a}. Then replacing the minimum over $\tilde{\Theta}$ by setting $\tilde{\theta} = \hat{\theta}$ and rearranging, one can conclude, for the estimator given by 
\begin{align*}
\hat{\theta}(X) = \argmin_{\theta \in \Theta} \left(-\log P_{\theta}(Z) +  pen(\theta)\right)
\end{align*}
redundancy risk bounds as follows:
\begin{align}{\label{risk}}
\E B_{tot}(P,P_{\hat{\theta}}) \leq \E \min_{\theta \in \Theta} \left(\log \frac{P(Z)}{P_{\theta}(Z)} + pen(\theta)\right).
\end{align}
Here $B_{tot}$ is the total Bhattacharya divergence(see \cite{Bhattacharyya1943}) and in i.i.d cases~\eqref{risk} becomes 
\begin{align*}
\E B(P,P_{\hat{\theta}}) \leq \frac{1}{n}\E \min_{\theta \in \Theta} \left(\log \frac{P(Z)}{P_{\theta}(Z)} + pen(\theta)\right),
\end{align*}
where $B$ is the single sample Bhattacharya divergence and $n$ is the sample size. When we say $\tilde{pen}(\tilde{\theta})$ is Kraft summable we mean
\begin{equation}{\label{Kraft}}
\sum_{\tilde{\theta} \in \tilde{\Theta}} \exp(-\L(\theta)) \leq 1. 
\end{equation}
\begin{remark}
The redundancy which is an expected minimum excess codelength, can be further upper bounded by the minimum expected excess codelength which is called the index of resolvability as in \cite{BarronCover1991}.
$$ Res = \min_{\theta \in \Theta} \left(D(P,P_{\theta}) + pen(\theta)\right).$$
Hence, we have the relation
\begin{equation}{\label{target}}
Risk \leq Redundancy \leq Resolvability.
\end{equation} 
As we see from the definition of Resolvability, the upper bound of the risk is governed by an ideal tradeoff between Kulback approximation error and complexity. So, in this sense the two stage estimator is adaptive, it looks at the tradeoff between approximation error and complexity at the population size relative to the sample size given. 
\end{remark}
We call the penalties for which the penalized likelihood procedure gives risk bounds of the form $$ Risk \leq Redundancy $$ as statistical risk valid penalties. It is also shown in \cite{BarronHuangLiLuo2008a} \cite{BarronHuangLiLuo2008b} that the $l_1$ penalty in linear regression and log density estimation problems is indeed codelength valid and statistical risk valid. We add a new example to this story in this paper where we present the risk validity of the $l_1$ type penalty in Gaussian Graphical Models. 
Also in this paper, we present a new interpretation of the traditional $l_0$ penalty in problems such as linear or logistic regression being codelength valid. Traditionally in the MDL  literature, the penalty $pen(\theta) = \frac{dim(\theta)}{2} \log(n) + o(n)$ in linear regression has been shown to be codelength valid but with the drawback that the $o(n)$ term is unbounded with respect to $\theta$ as we go out to the edges of the parameter space. Here we resolve this issue as we show that the penalty $2 \frac{dim(\theta)}{2} \log(n) + o(n)$ can be indeed thought of as codelength valid in the regime $n>p$ with the $o(n)$ term remaining bounded as a function of $\theta.$ Also in the linear regression problem, by slightly modifying the existing risk bound techniques, we are able to show the same penalty, described before found to be codelength valid, is also statistically risk valid.
\begin{remark} We comment that the requirement that the penalizer be $2$ times a codelength valid penalty in~\eqref{verifyrisk} stems from the loss function being the Bhattacharya divergence. The theory goes through for any of the Chernoff-Renyi divergences with parameter $\alpha$ strictly between $0$ and $1$, in which case the factor of $2$ may be replaced by 
$\frac{1}{\alpha}.$
\end{remark}

\section{Validity of $l_1$ penalty}
In order to show validity of a $l_1$ type penalty, the countable subset $\tilde{\Theta}$ is taken to be the $\delta$ integer lattice and now a penalty  is defined on $\tilde{\Theta}$
\begin{equation}{\label{pendef1}}
\tilde{pen}(\tilde{\theta}) = \frac{2}{\delta} \Vert \tilde{\theta} \Vert_1 \log(4p) + 2 \log 2.
\end{equation}
The actual value of $\delta$ depends on the problem at hand and is chosen accordingly. The fact that $\tilde{pen}(\tilde{\theta})$ defined this way satisfies~\eqref{Kraft} is by seeing it as $-2 \log$ of a subprobability measure, say $\pi$ on $\tilde{\Theta}.$ For any non negative integer $k,$ $\pi$ puts mass $(1/2)^{k+1}$ on the set of all points on the lattice with $l_1$ norm equal to $k$ times $\delta.$ Also a counting argument shows that there are at most $(2p)^k$ points in the $\delta$ lattice with $l_1$ norm equal to $k$ times $\delta.$ So, setting $\pi$ at any point in the grid to be equalling $\frac{1}{(2p)^k}$ times $(1/2)^{k+1},$ where $k$ is the $l_1$ norm of the particular point divided by $\delta,$ defines a subprobability. Now if one shows an equivalent version of~\eqref{verifycode}, that for any given $\theta \in \Theta$ and data $X,$ one has 
\begin{align}{\label{verifycode2}}
\log P_{\tilde{\theta}}(Z) - \log P_{\theta}(Z) + \tilde{pen}(\tilde{\theta}) \leq pen(\theta) 
\end{align}
then one shows codelength validity of $pen(\theta).$ Analogously if one shows an equivalent version of~\eqref{verifyrisk}, that for any given $\theta \in \Theta$ and data $X,$ 
such that $pen(\theta)$ satisfies 
\begin{align}{\label{verifyrisk2}}
\begin{split}
&\min_{\tilde{\theta} \in \tilde{\Theta}} [\log P_{\tilde{\theta}}(Z) - \log P_{\theta}(Z) \\&- B(\theta^{\star},\tilde{\theta}) + B(\theta^{\star},\theta) + \tilde{pen}(\tilde{\theta})] \leq pen(\theta) 
\end{split}
\end{align}
then we get the desired risk bound~\eqref{risk} of the penalized likelihood procedure. The validity of $l_1$ penalty in the case of linear models and log density estimation has already been shown by techniques as above in \cite{BarronHuangLiLuo2008a} and \cite{BarronHuangLiLuo2008b}. We present new results in the case of Gaussian Graphical Models.
\subsection{Gaussian Graphical Models}
Consider the problem of estimating the inverse covariance matrix of a multivariate gaussian random vector. Suppose we observe $Z = (x_1,\dots,x_n),$ each of which is a $p$ length vector drawn i.i.d from $N(0,\theta^{-1}).$ We denote the true inverse covariance matrix to be $\theta^{\star}.$ Let us denote the $- \log det $ function as $\phi.$ This $\phi$ is a convex function on the space of all $p \times p$ matrices with the convention that $\phi$ takes value $+ \infty$ on any matrix that is not positive definite. Inspecting the log likelihood of this model we have 
\begin{align*}
\frac{1}{n} \log P_{\theta}(Z) = \frac{p}{2} \log(2 \pi) +\frac{1}{2} Tr(S \theta) + \frac{\phi(\theta)}{2}
\end{align*}
Here, $Tr(S \theta)$ is the sum of diagonals of the matrix $S \theta$ and $S = \frac{1}{n} \sum_{i=1}^{n} x_i^T x_i.$ In this setting $\theta_{ij} = 0$ means that the $i$th and $j$th variables are conditionally independent given the others. We outline the proof of the fact that the penalty $\vert \theta \vert_1,$ which is just the sum of absolute values of all the entries of the inverse covariance matrix, is a statistical risk valid penalty. For this model, the Bhattacharya Divergence is 
$$B(\theta_1,\theta_2) = \frac{1}{2} [\phi(\theta_1) + \phi(\theta_2)] - \phi([\theta_1 + \theta_2]/2).$$ 
We assume that the truth $\theta^{\star}$ is sufficiently positive definite in the following way. We assume that for any matrix $\{\Delta: \Vert \Delta \Vert_{\infty} \leq \delta \}$ we have  
\begin{equation}{\label{ass}}
(\theta^{\star} + \Delta) \succ 0.
\end{equation}
We remark that this is our only assumption on the true inverse covariance and the $\delta$ in the assumption is the same $\delta$ used in constructing the countable set. The value of $\delta$ is specified later. Now we give an idea as to how we verify~\eqref{verifyrisk2} and establish risk validity of the $l_1$ penalty times a multiplier.
Our strategy is to upper bound the left side in~\eqref{verifyrisk2} which is a minimum over the entire $\delta$ lattice by a minimum over the $2^p$ vertices of the cube that $\theta$ lives in. We further upper bound this minimum by an expectation over these vertices with the random choice being unbiased for $\theta.$ Taylor expanding the log likelihoods and the Bhattacharya divergence terms upto the second order permits us to do careful reasoning and obtain that the l.h.s in\eqref{verifyrisk2} can be upper bounded by $$ n D_{max}(\Sigma^{\star})  \Vert \theta \Vert_1 \delta + \frac{2}{\delta} \Vert \theta \Vert_1 \log(4p) + 2 \log 2.$$ $D_{max}(\Sigma^{\star}) $ is the square of the maximum diagonal of  the true covariance matrix $\Sigma^{\star}.$ We need assumption~\eqref{ass} to ensure we stay in the region where $\phi$ is smooth and differentiable. 
Now by setting $\delta  = \sqrt{\frac{2 \log(4p^2)}{n}},$ we again see that penalty defined as below satisfies~\eqref{verifyrisk2}.
\begin{align*}
\begin{split}
&pen(\theta) = 2 D_{max}(\Sigma^{\star}) \sqrt{2n \log(4p^2)}) \Vert \theta \Vert_1 \\&+ 2 \log 2
\end{split}
\end{align*} 
So in this case a suitable multiplier times the $l_1$ penalty plus a constant is a risk valid penalty. We present the risk bound after dividing throughout by $n.$

\begin{theorem}
Define the estimator as
\begin{align*}
\begin{split}
\hat{\theta} = \argmin_{\theta \in \R^p} &[\frac{1}{2} Tr(S \theta) + \phi(\theta) \\&+ \lambda {\Vert \theta \Vert}_1]. 
\end{split}
\end{align*}
Then we have the risk bound
\begin{align*}
\begin{split}
&\E B(\theta^{\star},\hat{\theta}) \leq \E \inf_{\R^{p}} [\frac{1}{n}\log \frac{P_{\theta^{\star}}(Z)}{P_{\theta}(Z)} + \lambda \Vert \theta \Vert_1]\\&+ \frac{2 \log 2}{n}
\end{split}
\end{align*}
whenever $\lambda \geq  2 D_{max}(\Sigma^{\star}) \sqrt{\frac{\log(4p^2)}{n}}.$ 
\end{theorem}
\begin{remark}
By setting $\theta = \theta^{\star}$ in the right side of the bound, as long as $\theta^{\star}$ has finite $l_1$ norm, one has the standard risk bound $\sqrt{\frac{\log(4p^2)}{n}} \Vert \theta^{\star} \Vert_1.$ The main purpose of the risk bound is to demonstrate the adaptation properties of the $l_1$ penalized estimator and to demonstrate redundancy, a coding notion, as the upper bound to the statistical risk which has been championed in \cite{Grunwald2007}.
\end{remark}
\begin{remark}
The assumption~\eqref{ass} says that the true inverse covariance matrix $\theta^{\star}$ should be in the interior of the cone of positive definite matrix by a little margin. This assumption may be acceptable even in high dimensions as it does not prohibit collinearity.
\end{remark}

\begin{flushleft}
\section{Validity of $l_0$ penalty in Linear Regression}
\end{flushleft}
The log likelihood of the model is $$-\log P_{\theta}(Z) = \frac{1}{2 \sigma^2} \Vert y - X \theta \Vert_2^{2} + \frac{n}{2} \log 2\pi \sigma^2$$ where $Z = (y_{n \times 1},X_{n \times p})$ and $X$ is the design matrix and $\sigma^2 I $ is the known covariance matrix of $y.$ Assuming the sample size $n$ is larger than the number of explanatory variables $p,$ we divide the data into $Z_{in} = (y_{in},X_{in})$ consisting of $p$ samples and $Z_{f} = (y_{f},X_{f})$ consisting of $(n-p)$ samples. We do assume that the $l_2$ norms of the columns of $X$ are upper bounded by $n.$ Then for $pen(\theta|Z_{in})$ having the leading term $k(\theta) \log n,$ we show, for a Kraft valid codelength $\L(\tilde{\theta}|Z_{in})$ on $\tilde{\Theta}$ the following inequality analogous to~\eqref{verifycode}.
\begin{align}{\label{verifycode3}}
\begin{split}
&\min_{\theta \in \Theta} \{-\log P_{\theta}(Z) + pen(\theta|Z_{in})\} \geq \\&\min_{\tilde{\theta} \in \tilde{\Theta}} \{-\log P_{\tilde{\theta}}(Z_{f}) + 2 \L(\tilde{\theta}|Z_{in})\}
\end{split}
\end{align}
where now the right side of \eqref{verifycode3} gives a two stage codelength interpretation provided we treat it as codelengths on $Z_{f}$ conditional on $Z_{in}$ and hence the left side as a function on $Z_{f},$ being not less than the right side, also has a two stage conditional codelength interpretation. We now sketch how to show~\eqref{verifycode3}. Let us make some relevant definitions. 
For $S \subset [1:p],$ let $X_{in,S}$ denote the initial part of the design matrix with column indices in $S.$ Hence $X_{in,S}$ is a $p$ by $\vert S \vert$ matrix. For $\theta \in \R^p$ we denote the support of $\theta$ or the set of indices where $\theta$ is non zero by $S(\theta).$ We denote $\vert S(\theta)\vert$ by $k(\theta).$ Let $S^{\star}$ be the support of the true vector of coefficients $\theta^{\star}.$ For any given $S,$ denoting $\int_{\R^{\vert S \vert}} (\frac{P_{\phi}(Z_{in)}}{P_{\theta^{\star}}(Z_{in})})^{1/2} d\phi $ by $Int(Z_{in},S)$ one can check that this integral is 
\begin{align}{\label{int1}}
\begin{split}
&\exp[\frac{1}{4}\!\left(\Vert y_{in}\!-\!X_{in,S^{\star}}\!\theta^{\star} \Vert_2^2\!-\!\Vert y_{in}\!- O_{X_{in,S}}\!y_{in}\!\Vert_2^2\right)] \\& \sqrt{4\pi^{\vert S \vert} det(X_{in,S}^T X_{in,S})^{-1}}  
\end{split}
\end{align} 
where $O_{X_{in,S}}$ denotes the orthogonal projection matrix onto the column space of the matrix $X_{in,S}.$ In order to have a lower approximation for the integral $\int_{\R^{\vert S \vert}} P_{\phi}(Z_{in}) d\phi $ one can divide $\R^{\vert S \vert}$ into cubes of sidelength $\delta > 0$ with vertices of the cube being the $\delta$ integer lattice and for each cube choose the point where $P_{\phi}(Z_{in})$ is minimized within the cube. The value of $\delta$ is specified later. We define the set $C_{S} \subset \R^{\vert S \vert}$ to be the set of points obtained this way. Although $C_{S}$ depends on $\delta$ we do not explicitly write it for convenience. As we have defined, $C_{S} \subset \R^{\vert S \vert}$ but by imagining the other coordinates to be zero we can treat $C_{S} \subset \R^p.$   A property of the set of points $C_{S}$ is that the Reimann sum is upper bounded by the integral. Indeed
\begin{equation}{\label{int2}}
\sum_{\phi \in C_{S}} (\frac{P_{\phi}(Z_{in)}}{P_{\theta^{\star}}(Z_{in})})^{1/2} (\frac{\delta}{2})^{\vert S \vert} \leq Int(S). 
\end{equation}
We now define the countable set $$\tilde{\Theta} = \cup_{k = 0}^{p} \cup_{\{S:\vert S \vert = k\}} C_{S}.$$ Define $h(\tilde{\theta},Z_{in})$ as
\begin{align*}
\begin{split}
&h(\tilde{\theta},Z_{in}) = (\frac{1}{2})^{k(\tilde{\theta})+1} \frac{1}{{p \choose k(\tilde{\theta})}} \\&(\frac{\delta}{2})^{k(\tilde{\theta})} (\frac{1}{\sqrt{2}})^{k(\theta^{\star})} \sqrt{1/(4\pi)^{k(\tilde{\theta})}}\\& \sqrt{det(X_{in,S(\tilde{\theta})}^T X_{in,S(\tilde{\theta})})} (\frac{P_{\tilde{\theta}}(Z_{in})}{P_{\theta^{\star}}(Z_{in})})^{1/2}.  
\end{split}
\end{align*}
We also define $$ M(Z_{in}) = \sum_{\tilde{\theta} \in \tilde{\Theta}} h(\tilde{\theta},Z_{in}).$$
Now by first summing over $k,$ then $\{S:\vert S \vert = k\}$ and then $C_{S}$ and from~\eqref{int1} and ~\eqref{int2} we have 
\begin{align*}
\begin{split}
&M(Z_{in}) \leq \sum_{k = 0}^{p} \sum_{\{S:\vert S \vert = k\}} (\frac{1}{2})^{k+1}\frac{1}{{p \choose k}} (\frac{1}{\sqrt{2}})^{\frac{k(\theta^{\star})}{2}} \\&\exp[\frac{1}{4}\!\left(\Vert y_{in}\!-\!X_{in,S^{\star}}\!\theta^{\star} \Vert_2^2\!-\!\Vert y_{in}\!- O_{X_{in,S}}\!y_{in}\!\Vert_2^2\right)].
\end{split}
\end{align*}
Note that by properties of orthogonal projections,
\begin{align*}
\begin{split}
 &\Vert y_{in} -  X_{in,S^{\star}} \theta^{\star} \Vert_2^2 - \Vert y_{in} - O_{X_{in,S}} y_{in} \Vert_2^2 \leq \\& \Vert y_{in} -  X_{in,S^{\star}} \theta^{\star} \Vert_2^2 -\Vert y_{in} - O_{X_{in,S \cup S^{\star}}} y_{in} \Vert_2^2 = \\ & \Vert O_{X_{in,S \cup S^{\star}}} y_{in} - X_{in,S^{\star}} \theta^{\star} \Vert_2^2.
 \end{split}
 \end{align*} 
Also $\Vert O_{X_{in,S \cup S^{\star}}} y_{in} - X_{in,S^{\star}} \theta^{\star} \Vert_2^2$ is a $\chi^2$ random variable with degree of freedom atmost 
$ \vert S \vert + \vert S^{\star} \vert.$ Hence computing the mgf of a $\chi^2$ random variable at $1/4$, we have 
\begin{align*}{\label{mgf}}
\begin{split}
&\P_{in} \exp(\frac{1}{4} \Vert P_{X_{in,S \cup S^{\star}}} y_{in} - X_{in,S^{\star}} \theta^{\star} \Vert_2^2) \leq \\&(\sqrt{2})^{\vert S \vert + \vert S^{\star}\vert}.
\end{split}
\end{align*}
So we have 
\begin{align*}
\begin{split}
&\P_{in} M(Z_{in}) \leq \sum_{k = 0}^{p} \sum_{\{S:\vert S \vert = k\}}(\frac{1}{2})^{k+1}\frac{1}{{p \choose k}} \\&(\frac{1}{\sqrt{2}})^{k(\theta^{\star})} (\sqrt{2})^{k + k(\theta^{\star})}.
\end{split} 
\end{align*}
After cancellations we can simplify to obtain
\begin{equation*}
\P_{in} M(Z_{in}) \leq \sum_{k = 0}^{p} (\frac{1}{\sqrt{2}})^{k+1}.
\end{equation*}
By summing up the geometric series and taking expectation inside the $\log$ we get
\begin{equation}{\label{excess}}
\P_{in} \log M(Z_{in}) \leq \log(\frac{1}{2 - \sqrt{2}}).
\end{equation}
Finally, we define our conditional codelengths $\L(\tilde{\theta}|Z_{in}) = -\log h(\tilde{\theta},Z_{in}) + \log M(Z_{in}). $
Then, by definition, it is Kraft summable. Now given $\theta$ and $Z$, we seek to upper bound the minimum over $\tilde{\Theta}$ of the following expression 
\begin{align*}
\begin{split}
&\log \frac{P_{\theta}(Z)}{P_{\theta^{\star}}(Z)} - \log \frac{P_{\tilde{\theta}}(Z_{f})}{P_{\theta^{\star}}(Z_{f})} + 2 \L(\tilde{\theta},Z_{in}). 
\end{split} 
\end{align*}
The upper bound will then be a codelength valid penalty as seen from~\eqref{verifycode3}. Now our strategy to minimize the expression is to restrict to minimizing over $\tilde{\theta} \in C_{S(\theta)}$ which cannot decrease the overall minimum. In that case, except the $\log$ likelihood terms all the other terms remain fixed and hence do not have to be minimized. So then, using similar probabilistic techniques as in the $l_1$ case and then setting $\delta = \frac{1}{\sqrt{n}}$ we see that a codelength valid penalty $pen$ can be defined as
\begin{align*}
\begin{split}
&pen(\theta|Z_{in}) = k(\theta) \log(n) + 2 \log {p \choose k(\tilde{\theta})} \\& + \log \det X_{in,S(\theta)}^T X_{in,S(\theta)} + [1 + 4\log(2) - \\& \log(4\pi)] k(\theta) +\log(2) k(\theta^{\star}) + 2 \log M(Z_{in}).
\end{split} 
\end{align*}
The constant term $\log(2) k(\theta^{\star}) $ will not affect the optimization procedure. There are two terms which involve $Z_{in}.$ Now we assume that the maximum diagonal of $X_{in}^T X_{in}$ is at most $p.$ We think, this is reasonable in most situations.
Now using Hadamard's inequality we have
$\log det(X_{in,S(\theta)}^T X_{in,S(\theta)})$ term is at most $k(\theta) \log(p).$ Also as we have shown in~\eqref{excess}, the $2 \log M(Z_{in})$ term is only a constant on an average. So indeed, the dominating term of the penalty is $2$ times the traditional $k \frac{\log(n)}{2}$ term which we show has a conditional description length interpretation. 

We also have a risk bound arising from similar reasoning for a very similar penalty with leading term $k(\theta) \log n$ which is as follows:
\begin{align*}{\label{rbound}}
&\P B(P_{\theta^{\star}},P_{\hat{\theta}}) \leq (\frac{2}{n-p}) \log(\frac{1}{2 -\sqrt{2}}) \\&+ (\frac{1}{n-p}) \P \min_{\theta \in \Theta} \left(\log \frac{P_{\theta^{\star}}(Z)}{P_{\theta}(Z)} + pen(\theta)\right)
\end{align*}
where $\hat{\theta}$ is obtained from optimizing 
$$ \hat{\theta} = \argmin_{\theta \in \theta} \left(\frac{1}{2} \Vert y - X \theta \Vert_2^2 + pen(\theta)\right).$$
\begin{remark}
The mgf of a chi square, not existing at $\frac{1}{2}$ is the reason we suffer an extra factor of $2$ in the penalty. The argument can be modified to replace the factor of $2$ by anything strictly bigger than 1 with a larger additive constant in the penalty. The codelength interpretation and the risk bound is useful even in the regime $p/n$ is a constant with $n\!\rightarrow\!\infty.$ 
\end{remark}



\begin{thebibliography}{1,88}
\bibitem{BarronCover1991} A.R.\ Barron and T.M.\ Cover, ``Minimum
  complexity density estimation,'' \emph{IEEE Trans.\ Inform.\
    Theory}. Vol.37, No.4, pp.1034--1054.  1991.

\bibitem{BarronHuangLiLuo2008a} A.R.\ Barron, C.\ Huang, J.Q.\ Li, X.\
  Luo, ``The MDL principle, penalized likelihoods, and statistical
  risk,'' In \emph{Festschrift for Jorma Rissanen}. Tampere University
  Press, Tampere, Finland, 2008.

\bibitem{BarronHuangLiLuo2008b} A.R.\ Barron, C.\ Huang, J.Q.\ Li, X.\
  Luo, ``MDL, penalized likelihood and statistical risk,'' \emph{ IEEE
    Information Theory Workshop}.  Porto Portugal, May 4-9, 2008.

\bibitem{BarronRissanenYu1998} A.R.\ Barron, J.\ Rissanen, and B.\ Yu,
  ``The minimum description length principle in coding and modeling,''
  \emph{ IEEE Trans.\ Inform.\ Theory}. Vol.44, No.6, pp.2743--2760.
  1998.  Special Commemorative Issue: Information Theory: 1948-1998.

\bibitem{Bhattacharyya1943} A.\ Bhattacharyya, ``On a measure of
  divergence between two statistical populations defined by
  probability distributions,'' \emph{ Bull.\ Calcutta Math.\ Soc.}
  Vol.35, pp.99--109. 1943.

\bibitem{Grunwald2007} P. Gr\"unwald, \emph{The Minimum Description
    Length Principle}. Cambridge, MA, MIT Press. 2007.




\end{thebibliography}
\end{document}